\documentclass[11pt, a4paper, twoside]{amsart}
\usepackage{mathrsfs, tabularx}
\usepackage[all]{xypic}
\usepackage[active]{srcltx}

\usepackage {amsthm, amsmath, amssymb, amscd, float, latexsym, times, epic, eepic}

\theoremstyle{plain}
\newtheorem{theorem}{Theorem}[section]
\newtheorem{lemma}[theorem]{Lemma}
\newtheorem{proposition}[theorem]{Proposition}

\theoremstyle{definition}
\newtheorem{definition}[theorem]{Definition}

\newcommand {\Set}[1] {\mathbb{#1}}
\newcommand{\setR}[0]{\Set{R}}

\newcommand{\setC}[0]{\Set{C}}

\newcommand{\dom}[0]{{\mathscr{D}}}

\newcommand{\sJ}[0]{j}

\newcommand{\pd}[2]{\frac{\partial #1}{\partial #2}}

\newcommand{\od}[2]{\frac{d #1}{d #2}}

\newcommand{\pdd}[3]{\frac{\partial^2 #1}{\partial #2\,\partial #3}}

\newcommand{\vfield}[1]{{\mathfrak X}( #1)}
\newcommand{\varepsint}[0]{(-\varepsilon, \varepsilon)}

\newcommand{\slaz}[0]{\setminus \{0\}}

\title{$k$-parameter geodesic variations }

\keywords{semispray, Jacobi equation, geodesic variation, Jacobi tensor, shape operator, Riccati equation}
\author[Bucataru]{Ioan Bucataru}
\address{Ioan Bucataru, Faculty of Mathematics, Al.I.Cuza University
B-dul Carol 11, Iasi, 700506, Romania}

\author[Dahl]{Matias F. Dahl}
\address{
Matias F. Dahl, Institute of Mathematics, P.O.Box
1100, 02015 Helsinki University of Technology, Finland }

\date{\today}

\begin{document}
\begin{abstract}
  Suppose $S$ is a semispray on a manifold $M$. We know that the
  complete lift $S^c$ of $S$ is a semispray on $TM$ with the property
  that geodesics of $S^c$ correspond to Jacobi fields of $S$. In this
  note we generalize this result and show how geodesic variations of
  $k$-variables are related to geodesics of the $k$th iterated
  complete lift of $S$.  Moreover, for sprays (that is, homogeneous
  semi\-sprays) we show how geodesic variations of $(n-1)$-variables are
  related to a natural generalisation of Jacobi tensors.
\end{abstract}
\maketitle

%\newpage
%\tableofcontents
%\newpage
\section{Introduction}
\label{sec:intro}
Suppose $S$ is a \emph{semispray} on a smooth manifold $M$ of
dimension $n$. Then $S$ is a vector field on $TM\slaz$ (the manifold
on non-zero vectors), and a \emph{geodesic} is a curve $c\colon I\to
M$ such that $c'' = S\circ c'$. The motivation for studying semisprays
is that these provide a minimal mathematical structure for studying
curves on $M$ that solve systems of $2$nd order ODEs (ordinary
differential equations). In this way, semisprays provide a unified
setting for studying geodesics in Riemann, Finsler and
Lagrange geometries and for affine and non-linear connections. 
See for example, \cite{Shen2001, BucataruMiron:2007}.

As in Riemann geometry, one can study the variation of geodesics for a
semispray, and this leads to a \emph{Jacobi equation} that describes
the infinitesimal behaviour of a geodesic variation. 
More precisely, any geodesic variation induces a Jacobi field, and
conversely, any Jacobi field on a compact interval can be represented
by a geodesic variation \cite{BucataruMiron:2007, BucataruDahl:2008}.
The purpose of this paper is to study the analogous representation of
geodesic variations of multiple variables.
In Riemann geometry we know that geodesic variations of $n-1$
parameters are related to Jacobi tensors \cite{EBurgSullivan}. Here, a
\emph{Jacobi tensor} is a $(1,1)$-tensor along a geodesic that
satisfies an analogue to the usual Jacobi equation. For similar
results in Lorentz geometry, see \cite{EhrlichJung:1998} and
\cite{Larsen1996}.

The main results of this paper are Theorem
\ref{CharacterizationOfJacobiFields} and Theorem
\ref{prop:corrJacobiTensors}.  In Theorem
\ref{CharacterizationOfJacobiFields} we derive a Jacobi equation for
geodesic variations of $k\ge 1$ variables in the setting
semisprays. The advantage of this result is that it holds for an
arbitrary semispray. On the other hand, the disadvantage is that the
Jacobi field will be a curve in the $k$th iterated tangent bundle and
has $n 2^{k-1}$ components.  The proof relies on working with iterated
complete lifts \cite{Lewis:2001:GeomMax, 
  BucataruDahl:2008:conjugate}. In Theorem
\ref{prop:corrJacobiTensors} we specialise to sprays (that is,
homogeneous semisprays) and to manifolds of dimension $n\ge 2$. In
this setting we show that geodesic variations of $n-1$ variables
correspond to Jacobi tensors on $M$. This latter result can be seen as
a natural generalisation of the results in \cite{EBurgSullivan} and
\cite{EhrlichJung:1998} described above.  The advantage with Jacobi
tensors is that they only depend on $(n-1)^2$ components, and for
$k=n-1$ and $n\ge 2$ we always have $(n-1)^2 < n 2^{k-1}$.
As an application, we show in Proposition \ref{prop:coords} how
invertible Jacobi tensors can be used to construct global coordinates
around a geodesic.  

For invertible Jacobi tensors, there is a
close relation to \emph{tensor Riccati equations}.  We conclude the
paper with Section \ref{sec:ricEq}, which studies this correspondence
in more detail.
One motivation is that in Riemann geometry the Riccati equation
is satisfied by the \emph{shape operator} of hypersurfaces evolving
under the geodesic flow  \cite{
 KowalskiVanhecke:1986, Shen2001b}. 
In physics, a motivation for the Riccati equation is
that its trace correspond to the \emph{Raychaudhuri equation} 
used to study the expansion (and collapse) of a family of geodesics \cite{EBurgSullivan, 
EhrlichKim:1994,
JeriePrince:2000}.
A related equation is also the complexified Riccati equation which
describes the behaviour of amplitude for a propagating wave packet in
hyperbolic equations like the wave equation or Maxwell's equations
\cite{KKL:2001, Kachalov:2005, Dahl2007:Riccati}.
Lastly, for time-dependendent semisprays, one can define a generalisation
of a shape operator for a geodesic vector field, and moreover, show
that this shape operator satisfies a tensor Riccati equation
\cite{CrampinPrince:1984, JeriePrince:2000, JeriePrince:2002}.
In Proposition \ref{prop:RelationToJP} we show how this generalised shape operator
(for sprays) can be written explicitly using a geodesic variation.

\section{Preliminaries}
\label{sec:Notation}
We assume that $M$ is a smooth manifold without boundary and with
finite dimension $n\ge 1$. By smooth we mean that $M$ is a
topological Hausdorff space with countable base that is locally
homeomorphic to $\setR^n$, and transition maps are $C^\infty$-smooth.
All objects are assumed to be $C^\infty$-smooth on their domains.
 
For $r\ge 1$, let $T^rM=T\cdots TM$ be the $r$:th \emph{iterated
tangent bundle}, and for $r=0$ let $T^0M=M$. For example, when $r=2$
we obtain the second tangent bundle $TTM$ \cite{Besse:1978}, and in general
\begin{eqnarray*}
 T^{r+1}M &=& TT^{r}M,  \quad r\ge 0.
\end{eqnarray*}

For a tangent bundle $T^{r+1}M$ where $r\ge 0$, we denote
the canonical projection operator by $\pi_{r}\colon T^{r+1}M \to
T^rM$.  Occasionally we also write $\pi_{TTM\to M}$, $\pi_{TM\to M}, 
\ldots$ instead of $\pi_0\circ \pi_1$, $\pi_0, \ldots$.  
For $x\in T^rM$ where $r\ge 0$ let also $T_x^{r+1}M = \pi_r^{-1}(x)$.
%\begin{eqnarray*}
%T_x^{r+1}M &=& \{\xi \in T^{r+1}M : \pi_{r}(\xi)=x\}.
%\end{eqnarray*}

%%
We always use canonical local coordinates (induced by local
coordinates on $M$) for iterated tangent bundles.
If $x^i$ are local coordinates for $T^rM$ for some $r\ge 0$, we denote
induced local coordinates for $T^{r+1}M$, $T^{r+2}M$, and $T^{r+3}M$
by
\begin{eqnarray*}
%& & (x), \\
(x,y), \quad (x,y,X,Y), \quad (x,y,X,Y,u,v,U,V).
%& & (x,y,X,Y,u,v,U,V,a,b,A,B,p,q,P,Q).
\end{eqnarray*}
As above, we usually leave out indices for local coordinates and write
$(x,y)$ instead of $(x^i, y^i)$. 
%%
%In all the above local coordinates,
%$x^i$ are defined in a open set of $\setR^n$, and coordinates, $y^i,
%\ldots, Q^i$ are defined in $\setR^n$.

For $r\ge 1$, we treat  $T^rM$ as a vector bundle over the manifold
$T^{r-1}M$ with the vector space structure induced by projection
$\pi_{r-1}\colon T^rM\to T^{r-1}M$. Thus, if $\{ x^i : i=1, \ldots,
2^{r-1}n \}$ are local coordinates for $T^{r-1}M$, and $(x,y)$ are
local coordinates for $T^rM$, then vector addition and scalar
multiplication are given by
\begin{eqnarray}
\label{eq:TMplus}
(x,y) + (x, \tilde y) &=& (x,y+\tilde y), \\
\label{eq:mult}
\lambda \cdot (x,y)  &=& (x,\lambda  y).
\end{eqnarray}
For $r\ge 0$, a \emph{vector field} on an open set $U\subset T^rM$
is a smooth map $X\colon U\to T^{r+1}M$ such that $\pi_{r}\circ X =
\operatorname{id}_{U}$.  
The set of all such vector fields is denoted by $\vfield{U}$.

Suppose that $\gamma$ is a smooth map $\gamma\colon
\left(-\varepsilon, \varepsilon\right)^k \rightarrow T^rM $ where
$k\ge 1$ and $r\ge 0$. Suppose also that $\gamma(t^1,\ldots, t^k) =
(x^i(t^1, \ldots, t^k ))$ in local coordinates for $T^rM$.
Then the \emph{derivative} of $\gamma$ with respect to
variable $t^j$ is the curve $\partial_{t^j}\gamma \colon
\left(-\varepsilon, \varepsilon\right)^k$ $\to T^{r+1}M$ defined by
$\partial_{t^j}\gamma=\left(x^i, {\partial x^i}/{\partial
t^j}\right)$. When $k=1$ we also write $\gamma'=\partial_{t}\gamma$
and say that $\gamma'$ is the \emph{tangent of $\gamma$}.

Unless otherwise specified we always assume that $I$ is an open
interval of $\setR$ that contains $0$, and we do not exclude unbounded
intervals.  If $\phi\colon M\to N$ is a smooth map between
manifolds, we denote the tangent map $TM\to TN$ by $D\phi$, and if
%$\phi\colon M\to N$ is a smooth map between manifolds, 
$c\colon I\to M$ is a curve, then
\begin{eqnarray}
\label{curveDiff}
(\phi\circ c)'(t) = D\phi\circ c'(t), \quad t\in I.
\end{eqnarray}

\begin{lemma} 
\label{lemma:TrMRep}
If $\xi\in T^rM$ for some $r\ge 1$ then there exists a
 map 
\begin{eqnarray*}
  W\colon \varepsint^r&\to& M
\end{eqnarray*} 
such that if $s^1, \ldots, s^r$ are coordinates for $\varepsint^r$ then
\begin{eqnarray*}
   \xi = \partial_{s^1}\cdots \partial_{s^r} W(s^1, \ldots, s^r)\vert_{s^1, \ldots, s^r=0}.
\end{eqnarray*}
\end{lemma}

\begin{proof}
Let $V=\setR^{\dim M}$. For $k\ge 1$ let us define functions %$p_1, p_2, \ldots$ 
\begin{eqnarray*}
  w^{(k)}\colon V^{2^k}\times \setR^k \to V
\end{eqnarray*}
as follows.
For $k=1$ let
\begin{eqnarray*}
 w^{(1)}(u,v,s^1) &=& u+s^1 v, \quad u,v\in V, \,\, s^1\in \setR,
\end{eqnarray*}
and for $k\ge 2$, let 
\begin{eqnarray*}
 w^{(k)}(u,v,s^1, \ldots, s^k) &=& w^{(k-1)}(u+s^1 v, s^2, \ldots, s^k), \quad u,v\in V^{2^{k-1}},\\
& & \quad\quad\quad\quad\quad\quad \quad\quad \quad \quad\quad\quad\quad\!\!
 s^1, \ldots, s^k\in \setR.
\end{eqnarray*}
By induction, it follows that for all $k\ge 1$ we have
\begin{eqnarray*}
   \partial_{s^1} \cdots \partial_{s^k} w^{(k)}(u, s^1, \ldots, s^k)\vert_{s^1, \ldots, s^k=0} &=& u
\end{eqnarray*}
for all $u \in V^{2^k}$ and $s^1, \ldots, s^k\in \setR$.  If
$\phi\colon U\to V$ is a chart around $\pi_{T^rM\to M}(\xi)$, where
$V$ is as above, then the result follows using equation
\eqref{curveDiff} since the $r$th-fold tangent map $D^{(r)}\phi\colon
T^rU\to V^{2^r}$ is a coordinate chart around $\xi$.
\end{proof}

\subsection{Canonical involution on $T^rM$}
\label{CanonicalInvolution}
%\begin{definition}[Canonical involution on $T^rM$]
For $r\ge 2$, the \emph{canonical involution} is the unique
diffeomorphism $\kappa_r \colon T^rM\to T^rM$ that satisfies
\begin{eqnarray}
\label{KappaDefEq}
\partial_t \partial_s c(t,s) &=& \kappa_r\circ \partial_s \partial_t c(t,s)
\end{eqnarray}
for all curves $c\colon \varepsint^2 \to T^{r-2}M$. For $r=1$,
we define $\kappa_1=\operatorname{id}_{TM}$. For a discussion, see
\cite{BucataruDahl:2008} and references therein.
%\end{definition}

Let $r\ge 2$, let $x^i$ be local coordinates for
$T^{r-2}M$, and let $(x,y,X,Y)$ be local coordinates for $T^rM$. Then
%For a
%curve $c\colon \varepsint^2\to T^{r-2}M$, $c(t,s)=(x^i(t,s))$, we then
%have $\partial_t\partial_sc = \left(x^i, \pd{x^i}{s}, \pd{x^i}{t},
%\pdd{x^i}{t}{s} \right)$, and 
%In local coordinates, $\kappa_r$ is given by
\begin{eqnarray*}
%\label{localKappa}
\kappa_r(x,y,X,Y) &=& (x,X,y,Y).
\end{eqnarray*}
%In local coordinates for $TTM$ and $TTTM$ we have
%\begin{eqnarray*}
%  \kappa_2(x,y,X,Y) &=& (x,X,y,Y), \\
%  \kappa_3(x,y,X,Y,u,v,U,V) &=& (x,y,u,v,X,Y,U,V).
%%  \kappa_4(x,y,X,Y,u,v,U,V, a,b,A,B,p,r,P,R ) &=& \\
%  & & \!\!\!\!\!\!\!\!\!\!\!\!\!\!\!\!\!\!\!\!\!\!\!\!\!\!\!\!\!\!\!\!\!\!\!\!\!\!\!\!\!\!(x,y,X,Y,a,b,A,B, u,v,U,V,p,r,P,R).
%\end{eqnarray*}
For $r\ge 1$, we have identities
%\proofread{Equation \eqref{eq:DpiDDpi} follows by combining equations
%\eqref{commutationRelation} and \eqref{eq:pipikpipi} and taking the
%tangent map.}
\begin{eqnarray}
\label{eq:noNumberEq}
\kappa_r^2&=&\operatorname{id}_{T^rM},\\ % \quad r\ge 1, \\
\label{PiKappaXYZ}
\pi_r\circ D\kappa_r &=& \kappa_r\circ \pi_r, \\ %\quad r\ge 1, \\
\label{commutationRelation}
D\pi_{r-1} &=&  \pi_{r} \circ \kappa_{r+1}, \\ %\quad r\ge 2. \\
\label{piDDpi}
D\pi_{r-1}\circ \pi_{r+1} &=&  \pi_{r}\circ DD\pi_{r-1}, \\
\label{kappaId3}
DD\pi_{r-1} \circ \kappa_{r+2} &=& \kappa_{r+1} \circ DD\pi_{r-1}, \\
%\label{eq:pipikpipi}
%\pi_{r-1}\circ \pi_{r}\circ \kappa_{r+1} &=& \pi_{r-1}\circ \pi_{r}, \\
\label{eq:pipikpipiII}
\pi_{r-1}\circ D\pi_{r-1} &=& \pi_{r-1}\circ \pi_{r}.
\end{eqnarray}

%\subsection{Slashed tangent bundles $T^rM\slaz$}
%\label{sec:slashTrM}
The \emph{slashed tangent bundle} is the open set in $TM$ 
defined as
\begin{eqnarray*}
TM\slaz &=& \{ y\in TM : y\neq 0\}.
\end{eqnarray*}
On %iterated tangent bundles 
$T^rM$ for $r\ge 2$ we define 
\emph{slashed tangent bundles} as open sets
\begin{eqnarray*}
T^rM\slaz &=& \left\{ \xi \in T^rM : (D\pi_{T^{r-1}M\to M})(\xi)\in TM\slaz \right\}.
\end{eqnarray*}
For motivation, see Section \ref{sec:semispray}. Let also $T^rM\slaz = M$ when $r=0$. 

\subsection{Iterated lifts for functions}
\label{lifts}
Next we define the vertical and complete lift of a function $f\colon T^rM\to \setR$ on an iterated tangent bundle. When $r=0$, these lifts coincide
with the usual vertical and complete lifts defined in \cite{Yano1973}.

\begin{definition}%[Vertical lift] 
\label{def:vertLift}
For $r\ge 0$, the 
\emph{vertical lift} of a function $f\in$ $ C^\infty(T^rM\slaz)$ is the
 function $f^v\in C^\infty(T^{r+1}M\slaz)$ defined by
\begin{eqnarray*}
 f^v(\xi) &=&  f\circ \pi_{r} \circ \kappa_{r+1}(\xi), \quad \xi\in T^{r+1}M\slaz,
\end{eqnarray*}
and the \emph{complete lift} is the function $f^c\in C^\infty(T^{r+1}M\slaz)$
defined by
\begin{eqnarray*}
 f^c(\xi) &=& df\circ \kappa_{r+1}(\xi), \quad \xi\in T^{r+1}M\slaz.
\end{eqnarray*}
\end{definition}

Suppose $f\in C^\infty(T^rM\slaz)$ where 
 $r\ge 1$. If $x^i$ are local coordinates for $T^{r-1}M$, and
$(x,y,X,Y)$ are local coordinates for $T^{r+1}M$, then 
\begin{eqnarray*}
%\label{generalFV}
 f^v(x,y,X,Y) &=& f(x,X), \\
f^c(x,y,X,Y) &=& \frac{\partial f}{\partial x^a}(x,X)y^a + \frac{\partial f}{\partial y^a}(x,X) Y^a. \end{eqnarray*}

\section{Semisprays}
\label{sec:semispray}
The motivation for studying semisprays is that they provide a unified
framework for studying geodesics for Riemannian metrics, Finsler
metrics, non-linear connections, and Lagrange geometries. See
\cite{BucataruMiron:2007, Sakai1992, Shen2001}. Following
\cite{BucataruDahl:2008} we next define a
semispray on an iterated tangent bundle $T^rM$.

\begin{definition} %[Semispray] 
\label{def:spray}
Let $r\ge 0$.  A \emph{semispray} on $T^{r}M$ is a vector field $S\in
\vfield{T^{r+1}M\slaz}$ such that
$(D\pi_{r})(S)=\operatorname{id}_{T^{r+1}M\slaz}$.
%\begin{enumerate}
%\item 
%\label{it:SprayI}
%\end{enumerate}
\end{definition}

Let $S$ be a semispray $S\in \vfield{T^{r+1}M\slaz}$ for some $r\ge
0$.  If $(x,y,X,Y)$ are local coordinates for $T^{r+2}M$, then $S$ is
locally of the form
\begin{eqnarray}
\label{SprayDef}
  S(x,y) &=& \left(x^i,y^i,y^i,-2G^i\left(x,y\right)\right)\\
\nonumber
    &=& \left. y^i \pd{}{x^i}\right|_{(x,y)} - \left. 2 G^i(x,y) \pd{}{y^i}\right|_{(x,y)},
\end{eqnarray} 
where $G^i$ are functions $G^i\colon T^{r+1}U\slaz \to \setR$ for some open $U\subset M$.
%

%\subsection{Geodesics in $T^rM$}
%% Def: regular curve
Suppose $\gamma$ is a curve $\gamma\colon I \to T^rM$ where $r\ge
0$. Then we say that $\gamma$ is \emph{regular} if $\gamma'(t)\in
T^{r+1}M\slaz$ for all $t\in I$. When $r=0$, this coincides with the
usual definition of a regular curve, and when $r\ge 1$, curve $\gamma$ is
regular if and only if curve $\pi_{T^rM\to M}\circ \gamma \colon I\to M$ is
regular. 

\begin{definition} %[Geodesic]
\label{def:geodesic} 
Suppose $S$ is a semispray on $T^rM$ where $r\ge 0$. Then a regular curve
$\gamma \colon I\to T^rM$ is a \emph{geodesic} of $S$ if and only if
\begin{eqnarray}
\label{eq:geogeo}
  \gamma'' &=& S\circ \gamma'. %, \quad t\in I.
\end{eqnarray}
\end{definition}

Suppose $S$ is a semispray on $T^rM$ and locally $S$ is given by equation
\eqref{SprayDef}. Then a regular curve $\gamma\colon I\to T^{r}M$,
$\gamma=(x^i)$, is a geodesic for $S$ if and only if
\begin{eqnarray}
\label{eq:SGeo}
  \ddot x^i + 2 G^i\circ \gamma' &=& 0.
\end{eqnarray}

In Definition \ref{def:geodesic} we have defined geodesics on open
intervals.  If $\gamma$ is a curve on a closed interval we say that
$\gamma$ is a geodesic if $\gamma$ can be extended into a geodesic
defined on an open interval.

A semispray $S \in \vfield{T^{r+1}M\slaz}$ is a \emph{spray} if
$S$ further satisfies $[\setC_{r+1}, S] = S$, where $[\cdot, \cdot]$
is the Lie bracket and $\setC_{r+1}$ is the \emph{Liouville vector field}
$\setC_{r+1}\in \vfield{T^{r+1}M}$ defined as
\begin{eqnarray*}
\setC_r(\xi) &=& \partial_s ( (1+s) \xi)|_{s=0}, \quad \xi\in T^rM.
\end{eqnarray*}

Then Euler's theorem for homogeneous functions \cite{BCS:2000} implies
that functions $G^i$ are \emph{positively $2$-homogeneous}, that is,
if $(x,y)\in T^{r+1}M\slaz$, then
\begin{eqnarray*}
 G^i(x,\lambda y) &=& \lambda^2   G^i(x, y), \quad \lambda>0.
\end{eqnarray*}
Thus, if $\gamma$ is a geodesic for a spray $S$, the curve
$t\mapsto \gamma(A t + B)$ for constants $A>0$ and $B\in \setR$ is
again a geodesic (when defined). 
%This symmetry need not hold for a
%semispray.

\subsection{Complete lifts for a semispray}
\label{sec:completeLiftsSprays}
Suppose $S$ is a semispray on $M$.  As described in the introduction,
the complete lift of $S$ is a new semispray $S^{(1)}$ on $TM$.
The motivation for studying the complete lift is that geodesics of
$S^{(1)}$ are Jacobi fields of $S$, which in turn describe geodesic
variations of $S$.
%
%Thus a geodesic of $S^{(r)}$ will be a curve $J\colon I\to T^rM$ for
%any $r\ge 0$.  Also, we know that geodesics of $S^c$ are Jacobi fields
%of $S$, and these correspond to geodesic variations for $S$
%\cite{BucataruDahl:2008}.  Next we prove Theorem
%\ref{CharacterizationOfJacobiFields}, which describe how geodesics of
%$S^{(r)}$ and $r$-variable geodesic variations are related.
%
Next we define \emph{iterated complete lifts} for a semispray $S$ and
Theorem \ref{CharacterizationOfJacobiFields} will show how these are
related to geodesic variations of $k$ variables.
%on $M$ into a semisprays on $T^{r}M$ for any $r\ge 1$.  
%If
%$S$ is a semispray on $M$, we denote repeated complete lifts of $S$
%as follows:
%$$
 % S^{(0)} = S, \quad   S^{(1)} = S^c, \quad   S^{(2)} = S^{cc}, \quad S^{(3)} = S^{ccc}, \ldots.
%$$

The below definition for the complete lift can essentially be found in
\cite[Remark 5.3]{Lewis:2001:GeomMax}. For a further discussion about
related lifts, see \cite{BucataruDahl:2008}.

\begin{definition}%[Complete lift of semispray] 
\label{def:CompleteLiftOfSpray}
Suppose $S$ is a semispray on $M$. Then the
\emph{complete lifts} of $S$ are  semisprays $S^{(1)}$, $S^{(2)}$, $\ldots$
 on $TM, TTM, \ldots$ defined 
inductively as follows. For $r=0$, let $S^{(0)} = S$ and for $r\ge 0$, let
$S^{(r+1)}$ %\in \vfield{T^{r+2}M\slaz}$ 
be the semispray on $T^{r+1}M$ defined as
\begin{eqnarray*}
%\label{eq:Scdef}
  S^{(r+1)} &=& D\kappa_{r+2} \circ \kappa_{r+3} \circ DS^{(r)} \circ \kappa_{r+2}.
\end{eqnarray*}
%where $DS$ is the tangent map of $S$, 
%\begin{eqnarray*}
 % DS\colon T(T^{r+1}M\slaz) &\to& T^2(T^{r+1}M\slaz).
%\end{eqnarray*}
\end{definition}
By induction we see that all $S^{(0)}$, $S^{(1)}, \ldots $ are
semisprays. In fact, if $S^{(r)}$ is a semispray on $T^rM$ for some
$r\ge 0$, and if we write $S^{(r)}$ as in equation \eqref{SprayDef},
then %locally
\begin{eqnarray}
\nonumber S^{(r+1)} &=& \left(x,y,X,Y,X,Y,-2(G^i)^v, -2\left(G^i\right)^c \right) \\
      &=& \label{compSc}
 X^i \pd{}{x^i} + Y^i \pd{}{y^i}  -2 (G^i)^v \pd{}{X^i} -2 (G^i)^c
  \pd{}{Y^i},
\end{eqnarray}
whence $S^{(r+1)}$ is a semispray on $T^{r+1}M$.  

\begin{definition}
If $S$ is a semispray on $M$, then a \emph{Jacobi field} is a geodesic of $S^{(1)}$.
\end{definition}

A main motivation for the above definition will be given by Theorem \ref{CharacterizationOfJacobiFields} below.
Alternatively, from equation \eqref{compSc} we see that $S^{(1)}$ coincides with the usual
definition of the complete lift of a semispray on $M$ \cite{BucataruMiron:2007}.
Hence, the geodesic equation for $S^{(1)}$ coincide with the usual Jacobi equation
for $S$ in Riemann, Finsler, or Lagrange geometry \cite{BucataruDahl:2008}.

The \emph{geodesic flow} of a semispray $S^{(r)}$ on $T^rM$ for $r\ge
0$ is defined as the flow of $S^{(r)}$ as a vector field. The next
proposition shows how the geodesic flows of $S^{(0)}$, $S^{(1)}$,
$\ldots$ are related to each other \cite{BucataruDahl:2008}.  In
particular, if $S$ is complete (as a vector field), then all complete
lifts $S^{(1)}$, $S^{(2)}$, $\ldots$ are complete \cite{Yano1973}.
Let us also note that if $S$ is a spray, then all complete lifts are sprays.
%What is more, if $S$ is  
%metrizable (that is, $S$ is induced by a Lagrangian function) then all
%lifts are also metrizable \cite{BucataruDahl:2008}. \HOX{References}

\begin{proposition}%[Geodesic flow for iterated complete lifts]
\label{completeLiftFlow}
Suppose $S$ is a semispray on $M$ and $S^{(0)}, S^{(1)}, S^{(2)}, \ldots$ are as above.
Suppose furthermore that for each $r\ge 0$,
\begin{eqnarray*}
   \phi^{(r)} \colon \dom(S^{(r)} ) \to T^{r+1}M\slaz, %\quad r\ge 0
%\quad    \phi^c\colon \dom(S^c) \to T^{r+2}M\slaz, 
\end{eqnarray*}
is the geodesic flows of semispray $S^{(r)}$ with maximal domain
\begin{eqnarray*}
  \dom(S^{(r)})\subset T^{r+1}M \slaz \times \setR.
%\quad \dom(S^c)\subset T^{r+2}M \slaz \times \setR. 
\end{eqnarray*}
%\begin{enumerate}
%\item 
For all $r\ge 0$ we then have
\begin{eqnarray}
\label{eq:DSDSC}
%  \dom(S^c) &=& \left\{ (\xi,t)\in T^{r+2}M\times \setR : ((D\pi_{r})(\xi), t)\in \dom(S) \right\}, \\
\left( (D\pi_r)\times \operatorname{id}_\setR\right) \dom(S^{(r+1)}) &=& \dom(S^{(r)})
\end{eqnarray}
%\item 
and
\begin{eqnarray}
\label{eq:SCflow}
  \phi^{(r+1)}_t(\xi) &=& 
  \kappa_{r+2}\circ D\phi^{(r)}_t \circ \kappa_{r+2} (\xi),\quad (\xi,t) \in \dom(S^{(r+1)}),
\end{eqnarray}
where $D\phi^{(r)}_t$ is the tangent map of the map $\xi\mapsto \phi^{(r)}_t(\xi)$ for a fixed  $t$.
%\end{enumerate}
\end{proposition}

\subsection{$k$-parameter geodesic variations}
When $k=1$ the next definition reduces to the usual definition of a geodesic variation. 

\begin{definition}%[Geodesic variation]
\label{GeodesicVariation} 
Let $k\ge 1$ and let $c\colon I\to M$ be a geodesic for a semispray
$S$ on $M$.  Then a \emph{$k$-parameter geodesic variation of $c$} is
a map $V\colon I\times \varepsint^k \to M$ such that
\begin{enumerate}
%\item[\emi] The set $U$ is an open subset of $\setR^r$ that contains $0$. 
\item $V(t,0,\ldots, 0)=c(t)$ for all $t\in I$.
\item $t\mapsto V(t,s^1, \ldots, s^k)$ is a geodesic for all
  $(s^1,\ldots, s^k)\in \varepsint^k$.
\end{enumerate}
\end{definition}

\begin{theorem}
\label{CharacterizationOfJacobiFields}
Let $S$ be a semispray on $M$ and let $r\ge 1$.
\begin{enumerate}
\item 
\label{PropI}
If $V\colon I\times \varepsint^k \to M$ is a $k$-parameter geodesic
variation for some $k\ge 1$, then the curve $j\colon I\to T^rM$
defined as
\begin{eqnarray}
\label{DeltaVarExpI}
\sJ%(t)
&=& \partial_{s^{i_1}}\cdots \partial_{s^{i_r}}{V}|_{s^1,\ldots,s^k=0}% \quad t\in I 
%&=& \partial_{i^1}\cdots \partial_{i^r}{V}|_{s^1,\ldots,s^r=0}% \quad t\in I
\end{eqnarray}
is a geodesic of $S^{(r)}$. Here $s^1, \ldots, s^k$ are Cartesian coordinates for $\varepsint^k$
and $i_1, \ldots, i_r$ are indices for these coordinates
%coordinates for $\varepsint^k$ 
so that $1\le i_1, \ldots, i_r\le k$.

\item 
\label{PropII}
If $I$ is compact and $j \colon I\to T^rM$ is a geodesic of $S^{(r)}$,
 then $\sJ$ can be written as 
\begin{eqnarray}
\label{DeltaVarExpII}
\sJ%(t)
&=& \partial_{s^{1}}\cdots \partial_{s^{r}}{V}|_{s^1,\ldots,s^k=0}, \quad t\in I,
%&=& \partial_{i^1}\cdots \partial_{i^r}{V}|_{s^1,\ldots,s^r=0}% \quad t\in I
\end{eqnarray}
for some $r$-parameter geodesic variation $V\colon I^\ast \times \varepsint^r \to
 M$, where $I^\ast$ is an open subset with $I\subset I^\ast$.
\end{enumerate}
\end{theorem}

\begin{proof}
For part \ref{PropI}, let $j^{(1)}, \ldots, j^{(r)}$ be maps
\begin{eqnarray*}
  \sJ^{(p)} \colon I\times \varepsint^k  &\to& T^pM, \quad p=1,\ldots,r
\end{eqnarray*}
 defined as
%\begin{eqnarray*}
 %  \sJ^p &=& \partial_{s^{r-p+1}}\cdots  \partial_{s^r}{V}. 
%\end{eqnarray*}
$$ 
   \sJ^{(1)} = \partial_{s^{i_r}} V, \quad
   \sJ^{(2)} = \partial_{s^{i_{r-1}}} \partial_{s^{i_r}} V, \quad
\cdots, \quad 
   \sJ^{(r)} = \partial_{s^{i_1}} \cdots \partial_{s^{i_r}} V.
$$
By induction we next show that for all $p=1,\ldots, r$, 
\begin{eqnarray}
\label{InductionClaim}
   \partial_t^2 \sJ^{(p)} %(t,s^1, \ldots, s^p) 
&=& S^{(p)}\circ \partial_t \sJ^{(p)} %(t, s^1, \ldots, s^p), 
\quad \mbox{on } I\times \varepsint^k.
\end{eqnarray}
For $p=1$,
equations \eqref{curveDiff}, \eqref{KappaDefEq},
\eqref{eq:noNumberEq} and geodesic equation $\partial^2_t V=S(\partial_t V)$ yield
\begin{eqnarray*}
S^{(1)}(\partial_t \sJ^{(1)}) &=& D\kappa_2 \circ \kappa_3 \circ DS\circ \kappa_2 \circ \partial_t \partial_{s^{i_r}}  V \\
%                           &=&   D\kappa_2 \circ \kappa_3 \circ \partial_{s^{i_r}} (S(\partial_{t}  V)) \\
%&=&   \partial_t^2 \partial_{s^r} V \\
&=& \partial_t^2 \sJ^{(1)}.
\end{eqnarray*}
For $p\in \{1, \ldots, r-1\}$, the induction step follows by writing $j^{(p+1)}
= \partial_{s^{i_{r-p}}} j^{(p)}$ and repeating the above
calculation.
%but with different indices, 
%whence 
Part \ref{PropI} follows.
%$s^1,\cdots, s^r=0$ in equation \eqref{InductionClaim}.

For part \ref{PropII}, %suppose that $\sJ\colon I \to T^{r}M$ is a geodesic of $S^{(r)}$. By 
Lemma \ref{lemma:TrMRep} implies that
there exists a map $W\colon (-\delta, \delta)^{r+1} \to M$ with
\begin{eqnarray*}
   \sJ'(0) &=& \partial_{s^0} \cdots \partial_{s^{r}}  W|_{s^0, \ldots, s^{r}=0}
\end{eqnarray*}
With notation as in Proposition \ref{completeLiftFlow} we obtain
%For $k\ge 0$ let $\phi^{(k)}_t$ be the flow of $S^{(k)}$, and let
%$\dom(S^{(k)})\subset T^{k+1}M\slaz\times \setR$ be the domain of
%$\phi^{(k)}_t$. Then
\begin{eqnarray}
j(t) &=& \nonumber \pi_r\circ \phi^{(r)}_t\circ j'(0)\\
&=& \label{eq:jteq}
 \pi_r\circ \phi^{(r)}_t (\partial_{s^0} %\partial_{s^1}
 \cdots \partial_{s^{r}}  W|_{s^0, \ldots, s^{r}=0}), \quad t\in I.
\end{eqnarray}
We know that $\dom(S^{(r)})$ is open  \cite{AbrahamMarsden:1994}.
For each $t\in I$, we can therefore extend the domain of
$\phi^{(r)}_t \partial_{s^0} \cdots \partial_{s^{r}} W(s^0,\ldots,
s^r)$ to all $(t,s^0,\ldots, s^r)\in I_t\times (-\delta_t,
\delta_t)^{r+1}$ for some open interval $I_t\ni t$ and
$\delta_t>0$. Since $I$ is compact, we can extend $I$ into an open
interval $I^\ast\supset I$ and find an $\varepsilon^\ast>0$ such that
$\phi^{(r)}_t \partial_{s^0} \cdots \partial_{s^{r}} W$ is 
defined on $I^\ast \times (-\varepsilon^\ast,\varepsilon^\ast)^{r+1}$.
By equation \eqref{eq:DSDSC} it follows that for all $k\in \{0,\ldots, r\}$ we have
\begin{eqnarray*}
(\partial_{s^0} \partial_{s^{k+1}} \cdots \partial_{s^{r}} W,t)  &\in&\dom(S^{(r-k)}),\quad (t,s_1, \ldots, s_r) \in I^\ast \times
(-\varepsilon^\ast,\varepsilon^\ast)^{r}
\end{eqnarray*}
with convention $\partial_{s^{k+1}} \cdots \partial_{s^{r}} W = W$ for $k=r$.
For $k\in \{0,\ldots, r\}$, let 
\begin{eqnarray*}
j^{(k)}\colon   I^\ast \times (-\varepsilon^\ast,\varepsilon^\ast)^{r+1} &\to& T^rM
\end{eqnarray*}
be the map defined as 
\begin{eqnarray*}
j^{(k)}  &=&\partial_{s^1}\cdots \partial_{s^{k}} (\pi_{r-k} \circ \phi^{(r-k)}_t\circ \partial_{s^0} \partial_{s^{k+1}}\cdots \partial_{s^r} W(s^0, s^1, \ldots, s^r)).
%\vert_{s^1,  \ldots, s^r=0}, \quad t\in I.
\end{eqnarray*}
%For any $r\ge 1$,  
Equations \eqref{curveDiff}, \eqref{KappaDefEq},
\eqref{commutationRelation} and \eqref{eq:SCflow} imply $j^{(0)} =
\cdots = j^{(r)}$.  
Setting $s^0=\cdots =
s^r=0$ in equality $j^{(0)} = j^{(r)}$ and using equation
\eqref{eq:jteq} gives
\begin{eqnarray*}
%j(t)  &=&\partial_{s^1}\cdots \partial_{s^{r}} (\pi_0 \circ \phi^{(0)}_t\circ U)\vert_{s^1, \ldots, s^{r}=0}, \quad t\in I.
j(t)  &=& (\partial_{s^1}\cdots \partial_{s^{r}} V(t,s^1, \ldots, s^r))\vert_{s^1, \ldots, s^{r}=0}, \quad t\in I,
%\vert_{s^1, \ldots, s^r=0}, \quad t\in I.
\end{eqnarray*}
where 
   $V\colon I^\ast\times (-\varepsilon^\ast,\varepsilon^\ast)^{r}\to M$ 
is the geodesic variation
\begin{eqnarray*}
  V(t,s^1, \ldots, s^r) = \pi_0 \circ \phi^{(0)}_t\circ \partial_{s^0}  W(s^0,s^1, \ldots, s^r) \vert_{s^0=0}.
\end{eqnarray*}
Part \ref{PropII} follows.
 %suppose that $\sJ\colon I \to T^{r}M$ is a geodesic of $S^{(r)}$. By 
%\begin{eqnarray}
%\label{eq:RepeatApp}
 % \pi_r \circ \phi^{(r)}_t \circ \partial_a \partial_b S(a,b) = \partial_b(\pi_{r-1}\circ \phi^{(r-1)}_t \circ \partial_a S(a,b))
%\end{eqnarray}
%holds for any $S\colon (-\delta, \delta)^2 \to T^{r-1}M$ such that
%$\phi^{(r)}_t\partial_a\partial_b S(a,b)$ is defined.  Applying
%equation \eqref{eq:RepeatApp} $r$ times yields
\end{proof}

\subsection{Geodesics of $S^{(r)}$ and Jacobi fields of $S$}
\label{sec:vs}
The next two propositions describe how geodesics of $S^{(r)}$ for
$r\ge 2$ are related to geodesics of $S^{(0)}, \ldots, S^{(r-1)}$.  In
particular, Proposition \ref{ProjectionOfJacobiFields} shows that
geodesics of a semispray $S^{(r)}$ contain geodesics of all lower
order lifts $S^{(0)}, S^{(1)}, \ldots, S^{(r)}$. Thus, by equation
\eqref{eq:geogeo}, we can recover $S$ from any $S^{(r)}$ with $r\ge
1$.

\begin{proposition} 
\label{ProjectionOfJacobiFields}
Suppose $S$ is a semispray and $j\colon I\to T^rM$ is a geodesic of $S^{(r)}$.
\begin{enumerate}
\item 
\label{prop_ZaBi1}
If $r\ge 1$, then $\pi_{r-1} \circ j\colon I\to T^{r-1}M$ is a geodesic of $S^{(r-1)}$.
\item 
\label{prop_ZaBi2}
If $r\ge 2$, then $\kappa_r\circ j\colon I\to T^{r}M$ is a geodesic of $S^{(r)}$.
\item 
\label{prop_ZaBi3}
If $r\ge 2$, then $D\pi_{r-2} \circ j\colon I\to T^{r-1}M$ is  a geodesic of $S^{(r-1)}$.
\end{enumerate}
\end{proposition}

\begin{proof}
%  The claim in part \ref{prop_ZaBi1} is 
Since all claims are local  we may assume that
  $I$ is compact. Parts \ref{prop_ZaBi1} and \ref{prop_ZaBi2} follow by 
Proposition \ref{CharacterizationOfJacobiFields}. %  \ref{PropII}  and \ref{PropI}.
%
%there is a geodesic variation $V\colon I\times
%  \varepsint^r\to M$ such that
%\begin{eqnarray*}
%j(t)  &=& \partial_{s^1} \cdots \partial_{s^r} V(t,s^1, \ldots, s^r) |_{s^1, \ldots,s^r=0}.
%\end{eqnarray*}
%Hence
%\begin{eqnarray*}
%\pi_{r-1}\circ j(t)  &=& \partial_{s^2} \cdots \partial_{s^r} V(t,0, s^2, \ldots, s^r) |_{s^2, \ldots,s^r=0},
%\end{eqnarray*} 
%and part \ref{prop_ZaBi1} follows by Proposition \ref{CharacterizationOfJacobiFields} \ref{PropI}.
%%
%Part \ref{prop_ZaBi2} follows similarly, and 
Part \ref{prop_ZaBi3} follows by  equation \eqref{commutationRelation}.
%combining properties \ref{prop_ZaBi1} and
%\ref{prop_ZaBi2} and 
\end{proof}

The next proposition shows that every geodesic of $S^{(r)}$ induces
$r$ Jacobi fields for $S$. The proof of Proposition \ref{LastPart} 
follows by Proposition \ref{guessProjP} in Appendix \ref{sec:Canprojections}.

\begin{proposition}
\label{LastPart}
If $S$ is a semispray on a manifold $M$, and $j\colon I\to T^rM$ is a 
geodesic of $S^{(r)}$ where $r\ge 1$. Then there are $r$ distinct maps
\begin{eqnarray*}
 p^{(r)}_1, \ldots, p^{(r)}_r\colon T^rM &\to& TM
\end{eqnarray*} 
such that 
\begin{eqnarray*}
p^{(r)}_1\circ j, \ldots, p^{(r)}_r\circ j\colon I&\to& TM
\end{eqnarray*}
are geodesics of $S^{(1)}$ (that is, Jacobi fields of $S$).
\end{proposition}

Let us also note that if $S$ is a semispray on $M$ and $r\ge 0$, then
a geodesic $j\colon I\to T^{r}M$ for $S^{(r)}$ is uniquely determined
by $j'(t_0)\in T^{r+1}M\slaz$ for any $t_0\in I$. In contrast, Jacobi
fields $p^{(r)}_1\circ j, \ldots, p^{(r)}_r\circ j$ in the above
proposition do not determine $j$.  For example, suppose $j=(x,y,X,Y)$
is a geodesic of $S^{(2)}$, where $S$ is the flat spray on $M=\setR$.
Then $x,y,X,Y$ are independent and the geodesic equation reads
$$
    \ddot x = 0, \quad
    \ddot y = 0, \quad
    \ddot X = 0, \quad
    \ddot Y = 0.
$$
Now $p^{(2)}_1 \circ j = (x,y)$ and $p^{(2)}_2 \circ j = (x,X)$, but 
these do not determine $Y$.

\section{Jacobi tensors and geodesic variations}
\label{sec:jacobiTensors}
In the previous section, the main result was Theorem
\ref{CharacterizationOfJacobiFields}.  For a semispray, this theorem
shows how geodesics of the $k$th iterated complete lift are related to
geodesic variations of $k$-parameters.  Next we specialise this result
to sprays, that is, to semisprays that are homogeneous, and to
manifolds of dimensions $n\ge 2$.  The main result in this section is
Theorem \ref{prop:corrJacobiTensors}, which shows how geodesic
variations of $(n-1)$-parameters are related to \emph{Jacobi tensors}
(equation \eqref{eq:JacobiTensorEq}). Under additional assumptions,
this correspondence is known. For the Riemann case, see
\cite{EBurgSullivan} and for the Lorentz case, see
\cite{EhrlichJung:1998}. See also \cite{Larsen1996}.

Suppose $S$ is a semispray on a manifold $M$.  We know that $S$
induces a canonical \emph{dynamical covariant derivative} that
operates on arbitrary tensors on $TM\slaz$
\cite{BucataruDahl:Helmholtz}. We will here only need this derivative
for tensors along a geodesic \cite[Section 3.2]{BCD:2010} and for a
similar operator, see \cite[Definition 3.3]{JeriePrince:2002}.
Suppose $c\colon I\to M$ is a geodesic for $S$ and $X$ is a $(1,0)$-tensor along $c$.  If locally $X=X^i(t)\pd{}{x^i}\vert_{c(t)}$ we
define
\begin{eqnarray*}
  \nabla X &=& \left. \left( \od{X^i}{t} + N^i_j(c') X^j \right)\pd{}{x^i}\right\vert_{c(t)},
\end{eqnarray*}
where $N^i_j(y)=\pd{G^i}{y^j}(y)$ for $y\in TM\slaz$ and semispray $S$ is
written as in equation \eqref{SprayDef}.  Similarly, for a $(0, 1)$-tensor $\alpha=\alpha_i(t) dx^i\vert_{c(t)}$ we define
\begin{eqnarray*}
  \nabla \alpha &=& \left. \left( \od{\alpha_i}{t} - N^j_i(c') \alpha_j \right) dx^i\right\vert_{c(t)}.
\end{eqnarray*}
For a function $f\colon I\to M$ along $c$ we define $\nabla f =
\od{f}{t}$.  By the Leibnitz rule, the dynamical covariant  derivative $\nabla$ then extends to
tensors of any rank along $c$ \cite{BucataruDahl:Helmholtz}.  
For example, if $J$ is a $(1,1)$-tensor along $c$ and $v$ is a $(1,0)$-tensor, then
\begin{eqnarray}
\label{eq:nablaJ}
\nabla (J\circ v) &=& (\nabla J) \circ v + J\circ \nabla v.
\end{eqnarray}
We will say that a tensor $T$ along $c$ is \emph{parallel} if $\nabla
T = 0$.

For a semispray, the \emph{Jacobi endomorphism} is a $(1,1)$-tensor on
$TM\slaz$.  See \cite{BucataruDahl:Helmholtz} and references
therein. By restricting to $c'$ we define the \emph{Jacobi
  endomorphism} as the $(1, 1)$-tensor $\Phi$ along $c$ defined as
$\Phi(t) = \Phi^i_j (c') \pd{}{x^i}\otimes dx^j\vert_{c(t)}$, where
\begin{eqnarray*}
 \Phi^i_j &=& \left(2 \pd{G^i}{x^j} - S\left(\pd{G^i}{y^j}\right) - \pd{G^i}{y^r}\pd{G^r}{y^j}\right)_{c'(t)}.
\end{eqnarray*}
%Let us note that $1\choose 1$-tensors $J$ that satisfy equation \eqref{eq:JacobiT%ensorEq} are also called \emph{Jacobi tensors}.
If $S$ is a spray, then a geodesic $c\colon I\to M$ satisfies $\nabla c'=0$ and $\Phi(c')=0$.

\begin{definition} 
\label{def:jTensor}
Suppose $c\colon I\to M$ is a geodesic
for a semispray.
Then a \emph{Jacobi tensor}  %\colon \{c'\}^{\perp} \to \{c'\}^{\perp}$
  along $c$ is a $(1, 1)$-tensor $J$ along $c$
such that %$J\vert_{\{c'\}^\perp}$ 
\begin{eqnarray}
%\item $
\nabla^2 J + \Phi \circ J &=& \label{eq:JacobiTensorEq} 0.
%\item $
%(\nabla J)^t \cdot J - J^t\cdot (\nabla J)&=& \label{eq:LagEq} 0.
\end{eqnarray}
\end{definition}

As in the Riemannian case, Jacobi tensors and Jacobi fields are
related \cite{EBurgSullivan}: If $c\colon I\to M$ is a geodesic for a
semispray, a $(1,1)$-tensor $J$ along $c$ is a Jacobi tensor if and
only if $J\circ v\colon I\to TM$ is a Jacobi field for any parallel
vector field $v$ along $c$. See \cite[Proposition 2.10]{BCD:2010}.

When studying Jacobi tensors in the Riemann setting one usually
restrict them to tensors in the normal bundle $\{c'\}^\perp \to
\{c'\}^\perp$ and such Jacobi tensors can be characterised by their
initial values.  Our next goal is to prove Proposition
\ref{prop:parallelBasis}, which shows that a similar result
is true also for sprays. %To prove this result we first need some preliminary results.

Suppose $c\colon I\to M$ is a geodesic for a spray.  Suppose also that
$W$ is an $(n-1)$-dimensional subspace of $T_{c(0)}M$ and $c'(0)\notin
W$. 
We know that parallel transport is a linear isomorphism between tangent spaces.
Thus, by parallel transport we can extend
any basis $\{e_i \}_{i=1}^{n-1}$ for $W$ into linearly independent vectors
$\{e_i(t)\}_{i=1}^{n-1}$ in $T_{c(t)}M$ for any $t\in I$.  For $t\in
I$, let
\begin{eqnarray}
\label{eq:VtRep}
  W_t &=& \operatorname{span} \{e_1(t), \ldots,
e_{n-1}(t)\}.
\end{eqnarray}  Then $W_t$ does not depend on the choice of $\{e_i\}_{i=1}^{n-1}$, and
$T_{c(t)}M = \operatorname{span}\{c'(t)\}\oplus W_t$. Moreover, any 
vector field $v\colon I\to TM$ along $c$ can be written as
\begin{eqnarray}
\label{eq:jSmooth}
 v(t) &=& v^0(t) c'(t) + \sum_{i=1}^{n-1} v^i(t) e_i(t), \quad t\in I
\end{eqnarray}
for some smooth functions $v^0, v^1, \ldots, v^{n-1}\colon I\to \setR$.

Suppose $J$ is a $(1,1)$-tensor along a geodesic $c\colon I\to M$ for
a spray.  Moreover, suppose that $W\subset T_{c(0)}M$ is a subspace
such that $c'(0)\notin W$ and extensions $\{W_t: t\in I\}$ are defined
as above.  Then we say that $J$ is a \emph{transversal tensor with
  respect to $W_t$} if $J\circ c'(t)=0$ and $\operatorname{Im}
J(t)\subset W_t$ for all $t\in I$.
Equations \eqref{eq:nablaJ} and \eqref{eq:jSmooth} imply that $\nabla
J$ is transversal with respect to $W_t$ if $ J$ is transversal
with respect to $W_t$.

\begin{proposition} 
\label{prop:parallelBasis}
Suppose $c\colon I\to M$ is a geodesic for a spray and $J$ is a Jacobi
tensor along $c$.  Furthermore, suppose $W$ is an $(n-1)$-dimensional
subspace of $T_{c(0)}M$ with $c'(0)\notin W$.  
If 
\begin{enumerate}
\item $\Phi$ is tranversal with respect to $W_t$ and
\label{eq:transChar}
\item $J\circ c'(0) = 0, \,\,\, 
(\nabla J)\circ c'(0) = 0, \,\,\, 
\operatorname{Im} J(0)\subset W, \,\,\, 
\operatorname{Im} (\nabla J)(0)\subset W$ \!\!\!\!\!\!\!\!\!\!\!\!\!\!\!\!\!\!\!
\end{enumerate}
then $J$ is transversal with respect to $W_t$.
\end{proposition}

\begin{proof} 
%  Suppose $J$ is transversal.  If $v\colon I\to TM$ is a parallel
%  vector field along $c$, we can write $J\circ v$ as in equation
%  \eqref{eq:jSmooth} and $\operatorname{Im} (\nabla J)(0)\subset W$
%  follows by differentiation.  The other conditions in equation
%  \eqref{eq:transChar} are clear or follow similarly.
%
%  Conversely, 
Let
  $j\colon I\to TM$ be the vector field $j = J\circ c'$.  Then $j$ is
  a Jacobi field with $j(0)=\nabla j(0)=0$. Thus $j=0$.  To
  complete the proof we need to show that if $v\colon I\to TM$ is a
  parallel vector field along $c$, then $J\circ v\in W_t$ for all
  $t\in I$.  This follows by writing $j=J\circ v$ as in equation
  \eqref{eq:jSmooth} and computing the $c'(t)$-component of the Jacobi
  equation $\nabla^2 j + \Phi\circ j=0$.
\end{proof}

The next theorem shows how Jacobi tensors and geodesic variations are
related in the setting of sprays. For semisprays the analogous result
is Theorem \ref{CharacterizationOfJacobiFields}. When $n\ge 2$, 
the correspondence is not unique in neither theorem.

\begin{theorem}
\label{prop:corrJacobiTensors}
Suppose $S$ is a spray on a manifold $M$ of dimension $n\ge 2$.
\begin{enumerate}
\item 
\label{pp:paI}
Suppose $V\colon I\times \varepsint^{n-1}\to M$ is a geodesic
variation of $(n-1)$ parameters and $(t, s^1, \ldots, s^{n-1})$ are
coordinates for the domain. 
Furthermore, suppose $W$ is an $(n-1)$-dimensional vector space in $T_{c(0)}M$, 
$c'(0)\notin W$ and $e_i$ are as in equation \eqref{eq:VtRep}.
Then conditions 
\begin{eqnarray}
 J\circ c'(t) &=&  \label{eq:JtensorVariationI} 0,\\
J\circ e_a(t) &=&
\label{eq:JtensorVariationII}
 (\partial_{s^a}V)(t,0,\ldots, 0), \quad a\in \{1, \ldots, n-1\}, \quad t\in I.
%\operatorname{Flow}_{t}\circ % B\circ
%\pi\circ 
%P_{c(t)\to c(0)}(v), \quad v\in T_{c(t)} M,
\end{eqnarray}
define a  Jacobi tensor $J$ along $c$. Here $c(t)=V(t,0,\ldots, 0)$.
\item
\label{pp:paII}
Suppose $c\colon I\to M$ is a geodesic, where $I$
 is compact. Furthermore, suppose $J$ is a transversal Jacobi tensor with
respect to $W_t$ where $W_t$ is as in equation \eqref{eq:VtRep}
for some parallel vector fields $\{e_i(t)\}_{i=1}^{n-1}$.  Then there
exists an open interval $I^\ast \supset I$ and a geodesic variation
$V\colon I^\ast \times \varepsint^{n-1}\to M$ such that equations
\eqref{eq:JtensorVariationI}--\eqref{eq:JtensorVariationII} hold.
\end{enumerate}
\end{theorem}

\begin{proof}  
  For part \ref{pp:paI}, conditions
  \eqref{eq:JtensorVariationI}--\eqref{eq:JtensorVariationII} define a
  smooth $(1,1)$-tensor $J$ along $c$, and $J$ is a Jacobi tensor by Theorem
  \ref{CharacterizationOfJacobiFields} \ref{PropI}  and the observation after
  Definition \ref{def:jTensor}.

  For part \ref{pp:paII}, let $J_{a}\colon I\to TM$ be Jacobi fields $
  J_{a} = J(e_{a})$ for $a\in \{1, \ldots, n-1\}$. If $(x^i)_{i=1}^n$
  are local coordinates around $c(0)$ we can write $c(t)=(c^i(t))$ and
  $J_{a}(t) = J_{a}^i(t)\pd{}{x^i}\vert_{c(t)}$ for small $t$. In
  these coordinates, let $U\colon I\times \varepsint^{n-1}\to M$ be
  the map
\begin{eqnarray*}
  U(t,s^1, \ldots, s^{n-1}) &=& \left( c^i(0) 
    + t \dot{c}^i(0) 
    + \sum_{a=1}^{n-1} J_{a}^i(0) s^a 
    +\sum_{a=1}^{n-1} \dot{J}_{a}^i(0) ts^a\right)_{i=1}^n.
%\pi_0 \circ \phi^{(0)}_t\circ U(0,s^1, \ldots, s^{n-1})
\end{eqnarray*}
For all $ a\in \{1, \ldots, n-1\}$ it follows that
\begin{eqnarray*}
(\partial_t J_{a})(0) &=&   (\partial_t \partial_{s^a} U)(0,0,\ldots, 0).
\end{eqnarray*}
By a similar compactness argument as in the proof of Theorem
\ref{CharacterizationOfJacobiFields} \ref{PropII}, there exists an
open interval $I^\ast \supset I$ and an $\varepsilon>0$ such that the
map
\begin{eqnarray*}
 V(t,s^1, \ldots, s^{n-1}) &=&  \pi_0\circ \phi^{(0)}_t\circ  (\partial_t  U)(0,s^1, \ldots, s^{n-1}).
\end{eqnarray*}
defines a geodesic variation $V\colon I^\ast \times \varepsint^{n-1}\to M$.
Then 
\begin{eqnarray*}
  J_{a}(t) 
&=& \pi_1\circ \phi^{(1)}_t\circ  (\partial_t J_{a})(0) \\
&=& \pi_1\circ \phi^{(1)}_t\circ  (\partial_t \partial_{s^a} V)(0,0,\ldots, 0)\\
&=& \left(\partial_{s^a} V \right) (t,0,\ldots,0), \quad a\in \{1, \ldots, n-1\}, \quad t\in I
\end{eqnarray*}
and  part \ref{pp:paII} follows.
\end{proof} 

\subsection{Invertible Jacobi tensors}
\label{eq:Riccati}
Suppose $c$ is a geodesic $c\colon I\to M$ for a spray and $J$ is a
$(1,1)$-tensor along $c$.  Suppose also that $J$ is transversal with
respect to $W_t$ for some vector subspaces $\{W_t:t\in I\}$ as in
equation \eqref{eq:VtRep}.  Then we say that $J$ is \emph{invertible}
if $J\vert_{W_t}\colon W_t\to W_t$ is an invertible linear map for each
$t\in I$. By $J^{-1}$ we then denote the transversal
$(1,1)$-tensor determined by $(J^{-1})\vert_{W_t} =
(J\vert_{W_t})^{-1}$. 

For a semispray $S$ we say that points $p,q\in M$
are \emph{conjugate points} if there exists a Jacobi field $j\colon
I\to TM$ along a geodesic that connects $p$ and $q$, Jacobi field $j$
vanishes at $p$ and $q$, but $j$ is not identically zero.  Suppose $S$
has no conjugate points, and $J$ is a Jacobi tensor defined on
$I\subset \setR$ such that $J(0)=0$ and $\nabla J(0) =
\operatorname{Id}$.  Then $J$ is invertible on $I\slaz$. See
\cite{EhrlichKim:1994}.  
%%
%\HOX{This is clear: If $\det J(t_0)=0$ take a
%  $v \in Ker (J(t_0))$ and extend by parallel transport and consider
%  Jacobi field $j= J\circ v$ which has $j(0)=0$ and $j(t_0)=0$.}  
%%
In
the Riemann case, it also holds that Jacobi tensors that in addition
are \emph{Lagrange tensors} are invertible except at isolated points.
See \cite{EBurgSullivan}, and for a discussion of related results in
the Lorentz setting, see \cite{EhrlichKim:1994}.

Proposition \ref{prop:coords} below shows how the existence of an
invertible Jacobi tensor along a geodesic $c\colon I\to M$ implies
that geodesics near $c(I)$ can be straightened out and used to define
local coordinates around $c(I)$.  This gives sufficient conditions
when a geodesic variations is a diffeomorphism onto its range.  For a
similar result for vector fields, see \cite[Theorem
2.1.9]{AbrahamMarsden:1994}.

\begin{proposition} 
\label{prop:coords}
Suppose $c\colon I\to M$ is a geodesic for a spray $S$ on a manifold
$M$ with dimension $n\ge 2$, and suppose that $I$ is compact and
$c$ does not intersect
itself.
Furthermore, suppose $J$ is an invertible
Jacobi tensor along $c$.  Then there exists an open interval $I^\ast
\supset I$, an $\varepsilon>0$ and a map $V\colon
I^\ast\times \varepsint^{n-1}\to M$ such that
\begin{enumerate}
\item \label{prop:coords:I}
$V(t,0,\ldots, 0)=c(t)$ for $t\in I$, 
\item \label{prop:coords:II}
$V$ is a diffeomorphism onto its range, that is, $I^\ast\times
\varepsint^{n-1}$ are local coordinates around $c(I)$,

\item \label{prop:coords:III} for $(s^1, \ldots, s^{n-1})\in \varepsint^{n-1}$ the curve
\begin{eqnarray*}
t\mapsto V(t,s^1, \ldots, s^{n-1}), \quad t\in I^\ast
\end{eqnarray*}
is a geodesic.

\end{enumerate}
\end{proposition}

\begin{proof}
  Let $\{e_i(t)\}_{i=1}^{n-1}$ be a parallel basis such that $J$ is
  transversal with respect to
  $\operatorname{span}\{e_i(t)\}_{i=1}^{n-1}$.  Then Theorem
  \ref{prop:corrJacobiTensors} \ref{PropII} implies that there exists an open
  interval $I^{\ast}\supset I$, an $\varepsilon>0$ and a geodesic variation $V\colon
  I^{\ast} \times \varepsint^{n-1}\to M$ such that $V(t,0,\ldots,
  0)=c(t)$ for $t\in I$ and
\begin{eqnarray*}
\{ \partial_t V, \partial_{s^1} V, \ldots, \partial_{s^{n-1}} V\}
\end{eqnarray*}
are linearly independent for $t\in I$ and $s^1, \ldots, s^{n-1}=0$.  By
the inverse function theorem, we can restrict $I^\ast \supset I$ and
$\varepsilon>0$ such that $V$ is a local diffeomorphism. Since $c$
does not intersect itself, we can further apply \cite[Lemma
9.19]{spivakI} and restrict $I^\ast\supset I$ and $\varepsilon>0$ so
that $V$ is a bijection (and hence a diffeomorphism) onto its range.
\end{proof}

Suppose $S$ is a spray on $M$. Then $S$ is a \emph{Berwald spray} (or
\emph{affine spray}) if in coordinates $(x^i,y^i)$ for $TM\slaz$,
Christoffel symbols $\Gamma^i_{jk} = \pdd{G^i}{y^j}{y^k}$ do not
depend on $y^i$. This condition does not depend on coordinates
\cite[Section 6.1]{Shen2001}.

If $S$ is a Berwald spray in Proposition \ref{prop:coords}, then part
\ref{prop:coords:III} implies that $\Gamma^i_{11}=0$. That is, in the
terminology of \cite{MikesVanzurova:2010}, Proposition
\ref{prop:coords} gives sufficient conditions for the existence of
\emph{pre-semigeodesic coordinates} around a geodesic.
Let us note that a similar coordinate system is \emph{Fermi
  coordinates} around a geodesic, where all Christoffel symbols
satisfy $\Gamma^i_{jk}=0$, but only on the geodesic.  See
\cite{ManasseMisner:1963} and for the setting of Berwald sprays, see
\cite[p.~133]{Hicks:1965} and \cite[p.~64]{Eisenhart}.

\subsection{Riccati equation }
\label{sec:ricEq}
Suppose $J$ is a Jacobi tensor along a geodesic $c$ for a spray.
If  $J$ is invertible (and hence transversal), then
$L=\nabla J\circ J^{-1}$ is another transversal tensor along $c$ and
\begin{eqnarray}
\label{eq:riccati}
   \nabla L + L^2 + \Phi &=& 0.
\end{eqnarray}
This is the \emph{Riccati equation} for a $(1,1)$-tensor, and the
above observation demonstrates the relation between the Jacobi tensor
equation \eqref{eq:JacobiTensorEq} and the \emph{tensor Riccati
  equation} (equation \eqref{eq:riccati}).  In the Riemann (and
Finsler) case, the \emph{shape operator} of metric spheres satisfy the
Riccati equation \cite[Lemma 14.4.2]{Shen2001b}. For the Lorentz
setting, see \cite{EhrlichJung:1998}.  The trace of equation
\eqref{eq:riccati} is a generalisation of the \emph{Raychaudhuri
  equation} \cite{EhrlichKim:1994, JeriePrince:2000, JeriePrince:2002}.

Suppose $S$ is a semispray on a manifold $M$. Then the \emph{connection map} is 
the map $K\colon T(TM\slaz) \to TM$ defined as
\begin{eqnarray*}
  K(x,y,X,Y) &=& (x^i, Y^i + N^i_j(y) X^j),
\end{eqnarray*}
where $S$ is defined as in equation \eqref{SprayDef} and $N^i_j = \pd{G^i}{y^j}$.

\begin{definition}
Suppose $Z$ is a nowhere zero vector field $Z\in \vfield{U}$ defined on an open set $U\subset M$. % and suppose that $Z$ is nowhere zero.
Then $Z$ defines  a $(1,1)$-tensor on $U$ defined as
\begin{eqnarray*}
  A_Z &=& K \circ DZ. % (V), \quad V\in TU.
\end{eqnarray*}
\end{definition}
If $Z=Z^i\pd{}{x^i}$ in local coordinates $(x^i)$ we have
\begin{eqnarray*}
  A_Z &=& \left( \pd{Z^i}{x^j} + N^i_j(Z)\right) \pd{}{x^i}\otimes dx^j.
\end{eqnarray*}
From this local expression we see that tensor $A_Z$ in the above
corresponds to tensor $A_Z$ defined in Definition 4.7 in
\cite{JeriePrince:2000}.  
The next proposition is essentially \cite[Theorem 6.2]{JeriePrince:2000} formulated for geodesic variations.

\begin{proposition}
\label{prop:GeodRicc} 
Suppose $c\colon I\to M$ is a geodesic for a semispray 
on a manifold of dimension $n\ge 2$, 
suppose $V\colon I\times \varepsint^{n-1}\to M$ is a geodesic
variation of $c$ and $(t,s^1, \ldots, s^{n-1})$ are coordinates for the domain of $V$.
If $V$ is a diffeomorphism onto its range, then
\begin{eqnarray*}
   \nabla A_Z + A_Z^2 + \Phi &=& 0,
\end{eqnarray*} 
where $Z$ is the vector field induced by $\pd{}{t}$ and $A_Z$ is restricted to geodesic $c$.
\end{proposition}

\begin{proof}
The result follows by a direct computation and using that $\pd{G^i}{x^j}(c')=0$ in coordinates
$\{x^i\}_{i=0}^{n-1}$, where $x^0 = t$ and $x^i=s^i$ for $i=1, \ldots, n-1$.
\end{proof}

For sprays and semisprays we know that $A_Z$ acts as a generalisation
of the shape operator \cite{JeriePrince:2000, JeriePrince:2002} when
$Z=\partial_t V$ for a geodesic variation $V$.  In the setting of
sprays, the next proposition gives an explicit expression for $A_Z$.
In particular, the proposition gives sufficient conditions that imply
that $A_Z = \nabla J \circ J^{-1}$ for a suitable Jacobi tensor.  Let
us note that in Riemann and Finsler geometry, the shape operator for
small geodesic spheres can be written as $\nabla J\circ J^{-1}$ where
$J$ is a Jacobi tensor with $J(0)=0$ and $\nabla J =
\operatorname{Id}$.  See for example \cite{EBurgSullivan,
 KowalskiVanhecke:1986, Shen2001b}.  The next proposition gives
sufficient conditions when $A_Z$ has an analogous representation.  For
a discussion of equation \eqref{eq:J2}, see \cite{JeriePrince:2000,
  EBurgSullivan}.

\begin{proposition}
\label{prop:RelationToJP} 
%Suppose $S$, $V$ and $(t, s^1, \ldots, s^{n-1})$ are as in Proposition \ref{prop:GeodRicc}.
Suppose $c\colon I\to M$ is a geodesic for a spray  on a manifold of dimension $n\ge 2$, 
suppose $V\colon I\times \varepsint^{n-1}\to M$ is a geodesic variation of $c$,
suppose $(t,s^1, \ldots, s^{n-1})$ are coordinates for the domain of $V$,
and suppose  $W$ is an $(n-1)$-dimensional subspace in $T_{c(0)}M$ with $c'(0)\notin W$, 
and $J$ is a
$(1,1)$-tensor along $c$ defined as in Theorem \ref{prop:corrJacobiTensors} \emph{\ref{pp:paI}}.
Furthermore, suppose 
\begin{enumerate}
%$\operatorname{span}\{e^i\colon I_0\to TM\}_{i=1}^{n-1}$.
\item \label{t:II}
for some open interval $I_0\subset I$, the restriction $V\colon I_0\times \varepsint^{n-1}\to M$ is a diffeomorphism onto
its range. 
\item \label{t:I}
on $I_0$, $J$ is transversal with respect to $W_t$.
%$V$ is a diffeomorphism onto its range.
\end{enumerate}
Then $J$ is invertible on $I_0$, and  on $I_0$ we have
\begin{eqnarray}
A_Z &=&  \label{eq:J1} \nabla J\circ J^{-1},\\
%\nabla J &=& \label{eq:J1} A_Z\circ J, \\
\od{}{t} (\det J) &=& \label{eq:J2} \operatorname{trace} A_Z\cdot \det J,
\end{eqnarray} 
where $A_Z$ is the $(1,1)$-tensor along $c$ associated to vector field $Z
= \partial_t V$
and $det J$ is the determinant of the transverse part of $J$.
In particular, $A_Z$ is transversal on $I_0$.
\end{proposition}

\begin{proof}
  Let $e_1, \ldots, e_{n-1}$ be parallel vector fields along $c$ such
  that $J$ is defined by equations
  \eqref{eq:JtensorVariationI}--\eqref{eq:JtensorVariationII} and
  $W=\operatorname{span}\{e_a(0)\}_{a=1}^{n-1}$.  Then $J$ is
  transversal with respect to $\operatorname{span}\{e_a(t)\colon I_0\to
  TM\}_{a=1}^{n-1}$, and assumption \ref{t:II} implies that $J$ is
  invertible on $I_0$.  Let
\begin{eqnarray*}
    B &=& \{ J\circ v \colon I_0 \to TM : v\,\,\mbox{ is a parallel $(1,0)$-tensor along $c$ }\}.
\end{eqnarray*}
Then $B$ is a vector space over $\setR$ and $\operatorname{dim} B =
n-1$.  For $a\in \{1, \ldots, n-1\}$ we also have
$(\partial_{s^a}V)(t,0,\ldots, 0)\in B$.  Now $(t,s^1, \ldots,
s^{n-1})$ are local coordinates around $c(I_0)$, and any $j\in B$ can be written as
\begin{eqnarray*} 
 j(t) &=& \sum_{a=1}^{n-1} J^a \left. \pd{}{s^a}\right\vert_{c(t)}, \quad t\in I_0
\end{eqnarray*}  
for some constants $J^1, \ldots, J^{n-1}\in \setR$.
A direct computation shows that $\nabla j = A_Z\circ j$ on $I_0$ 
and equation \eqref{eq:J1} follows  by equation \eqref{eq:nablaJ}.
Equation \eqref{eq:J2} follows by the argument in \cite[Lemma 1]{EBurgSullivan},
that is, essentially by \emph{Liouville's formula}.
\end{proof}

\appendix
\section{Canonical projections $T^rM\to TM$}
\label{sec:Canprojections}
It is well known that there are two distinct canonical projections $TTM\to
TM$, 
namely $\pi_1, D\pi_0\colon TTM\to TM$. See for example \cite{Besse:1978, BucataruDahl:2008}. 
Similarly, for $TTTM$ there are three distinct canonical projections $TTTM\to TM$ 
illustrated in the below commutative diagram:
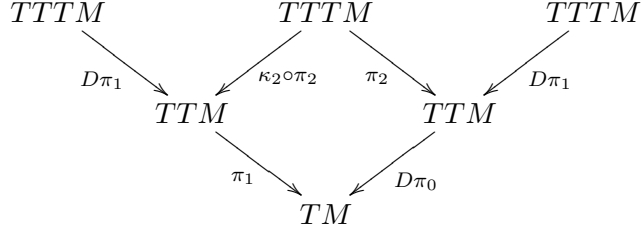
\begin{eqnarray*}
\begin{xy}
\xymatrix@C=1.0pc @R=2.0pc{
%\xymatrix{
TTTM & & TTTM & & TTTM \\
&  TTM \ar@{<-}[ul]^{D\pi_{1}}\ar@{<-}[ur]_{\kappa_2\circ \pi_{2}}
&  & TTM\ar[dl]^{D\pi_{0}} \ar@{<-}[ul]^{\pi_{2}}\ar@{<-}[ur]_{D\pi_{1}} \\
&  &  TM \ar@{<-}[ul]^{\pi_{1}}     &  & &}
\end{xy}
\end{eqnarray*}
Next we generalize these results and show that for any $r\ge 2$ there
are $r$ distinct canonical projections
$$
  p^{(r)}_1,\ldots,   p^{(r)}_r\colon T^rM\to TM.
$$ 
In this appendix we construct these projections are prove a number of 
technical properties. For geometric implications of these results, see Section \ref{sec:vs}.
%What is more, we show that $p^{(r)}_j$ are geodesic maps; that is,
%they map geodesics of $S^{(r)}$ into geodesics of $S^c$.  As a
%consequence, any geodesic $j\colon I\to T^rM$ of $S^{(r)}$ induces
%$r\ge 1$ Jacobi fields.

The maps are defined by induction. Let $p^{(1)}_1\colon TM\to TM$ be the identity map
\begin{eqnarray*}
  p_1^{(1)} &=& \operatorname{id}_{TM}.
\end{eqnarray*}
For $r\ge 2$ we define maps $p_1^{(r)},\ldots, p_k^{(r)}\colon T^rM\to TM$ as
\begin{eqnarray*}
p_i^{(r)} &=& 
\begin{cases}
p^{(r-1)}_i\circ \pi_{r-1}, & \mbox{for $i=1,\ldots, r-1$}, \\
p^{(r-1)}_{r-1}\circ D\pi_{r-2}, & \mbox{for $i=r$}.
%  D(\pi_{T^{k-1}M\to M}), & \mbox{for $i=k$}.
\end{cases}
\end{eqnarray*}
For example, for $TM, TTM, TTTM$ we obtain projection maps
\begin{eqnarray*}
p^{(1)}_1 &=& \operatorname{id}_{TM}, \\
p^{(2)}_1 &=& \pi_1, \quad   p^{(2)}_2 = D\pi_0, \\
p^{(3)}_1 &=& \pi_1\circ \pi_2, \quad   p^{(3)}_2 = D\pi_0\circ \pi_2, \quad   p^{(3)}_3 = D(\pi_{T^2M\to M}).
\end{eqnarray*}

\begin{proposition} 
\label{guessProjP}
Let $r\ge 1$.
\begin{enumerate}

\item 
\label{ProjMaps}
The maps  $p^{(r)}_1, \ldots, p^{(r)}_r \colon T^{r} M \to TM$
are distinct.

\item 
\label{projProp} 
%\begin{eqnarray}
$\pi_0\circ p_i^{(r)} = \pi_{T^rM\to M}$ for all $i=1, \ldots, r$.
%\end{eqnarray}

\item 
\label{projPropA} 
$p^{(r)}_1 = \pi_{T^rM\to TM}$ and $p^{(r)}_r = D(\pi_{T^{r-1}M\to M})$ for $r\ge 2$.
\item
\label{VariationPmap} 
Suppose $V\colon I\times \varepsint^k\to M$ is a map where $k\ge 1$. Then
\begin{eqnarray*}
    p^{(r)}_a\circ \partial_{s^1} \cdots \partial_{s^r} V%, |_{s^1, \ldots, s^r=0}
%(t,s^1, \ldots, s^k) 
&=& \partial_{s^{r-a+1}} V, %|_{s^1, \ldots, s^r=0}, %(t,s^1, \ldots, s^{n-1}), 
\quad a=1, \ldots, r,
\end{eqnarray*}
where $s^1, \ldots, s^r$ index coordinates for $\varepsint^k$ so that $1\le s^1, \ldots, s^r\le k$.
\end{enumerate}
\end{proposition}

\begin{proof} 
  For \ref{ProjMaps} the claim is clear for $r=1,2,3$ and for the
  induction step, suppose that \ref{ProjMaps} holds for some $r\ge 3$.
  Then
\begin{eqnarray}
\label{step_ABx2} 
p^{(r)}_i &\neq& p^{(r)}_{j},\ \ \mbox{ for all}\ \  1\le i<j\le r-1, \\
\label{step_ABx1} 
p^{(r)}_i &\neq& p^{(r)}_{r},\ \ \mbox{ for all}\ \  1\le i\le r-1.
\end{eqnarray}
Applying $\pi_r$ to equation \eqref{step_ABx2} from the right gives that 
$\{p^{(r+1)}_i\}_{i=1}^{r-1}$ are distinct. 
Similarly, applying $\pi_{r}$ to equation \eqref{step_ABx1} gives that 
$p^{(r+1)}_{i} \neq p^{(r+1)}_{r}$ for all $i=1, \ldots, r-1$,
Applying $D\pi_{r-1}$ to equation \eqref{step_ABx1} and using equation \eqref{eq:pipikpipiII} gives that 
$p^{(r+1)}_{i} \neq p^{(r+1)}_{r+1}$ for all $i=1, \ldots, r-1$.
Applying $DD\pi_{r-2}$ to equation \eqref{step_ABx1} with $i=r-1$ yields
\begin{eqnarray*}
   p^{(r-1)}_{r-1}\circ \pi_{r-1}\circ DD\pi_{r-2} &\neq&    p^{(r-1)}_{r-1}\circ D\pi_{r-2}\circ DD\pi_{r-2}
\end{eqnarray*}
whence equations \eqref{piDDpi} and \eqref{eq:pipikpipiII}
%to the left hand side and equation
%$D\pi_{r-1} \circ DD\pi_{r-1} = D\pi_{r-1}\circ D\pi_{r}$
%\eqref{eq:DpiDDpi} 
%to the right hand side 
imply that $p^{(r+1)}_{r} \neq
p^{(r+1)}_{r+1}$.
%%
%%
%% III
%%
Part \ref{projProp} follows by induction and equation  \eqref{eq:pipikpipiII}.
Part \ref{projPropA} follows directly by induction.
For part \ref{VariationPmap}, let us fix $k\ge 1$. The claim is clear for $r=1,2$
and for the induction step suppose that the claim holds for some $r\ge 2$.
For any $1\le s^1, \ldots, s^{r+1}\le k$ and $a=1,\ldots, r+1$ we then have
\begin{eqnarray*}
  p^{(r+1)}_a \partial_{s^1}\ldots \partial_{s^{r+1}}V
&=&  \begin{cases}
p^{(r)}_a\circ \pi_{r} \circ \partial_{s^1}\cdots \partial_{s^{r+1}}V, & a=1,\ldots, r, \\
p^{(r)}_{r}\circ D\pi_{r-1}\circ \partial_{s^1}\cdots \partial_{s^{r+1}}V, & a=r+1
\end{cases}\\
&=&  \begin{cases}
p^{(r)}_a\circ \partial_{s^2}\cdots \partial_{s^{r+1}}V, & a=1,\ldots, r, \\
p^{(r)}_{r}\circ \partial_{s^1}\partial_{s^3}\cdots \partial_{s^{r+1}}V, & a=r+1.
\end{cases}
\end{eqnarray*}
Part \ref{VariationPmap} follows by using the induction assumption in
the upper branch and by using %Proposition \ref{guessProjP}
part \ref{projPropA} in the lower branch.
\end{proof}

\subsection*{Acknowledgements}
We would like to thank Ricardo Gallego Torrom\'e for useful comments on the manuscript.

The work of IB was supported by a grant of the Romanian National
Authority for Scientific Research, CNCS UEFISCDI, project number
PN-II-RU-TE-2011-3-0017.  The work on MD has been supported by Academy
of Finland (project 13132527 and Centre of Excellence in Inverse
Problems Research) and by the Institute of Mathematics at Aalto
University.

%\cleardoublepage
%\bibliographystyle{abbrv}    %{amsalpha}
%\bibliographystyle{amsalpha}
%\bibliographystyle{apa}

%\nocite{*}    % list all entries
%\bibliography{../../../references}

\begin{thebibliography}{BCD11}

\bibitem[AM78]{AbrahamMarsden:1994}
R.~Abraham and J.E. Marsden, \emph{Foundations of mechanics}, Perseus books,
  1978.

\bibitem[BCD11]{BCD:2010}
I.~Bucataru, O.~Constantinescu, and M.F. Dahl, \emph{A geometric setting for
  systems of ordinary differential equations}, International Journal of
  Geometric Methods in Modern Physics \textbf{8} (2011), no.~6, 1291--1327.

\bibitem[BCS00]{BCS:2000}
D.~Bao, S.-S. Chern, and Z.~Shen, \emph{An introduction to {Riemann-Finsler}
  geometry}, Springer, 2000.

\bibitem[BD09]{BucataruDahl:Helmholtz}
I.~Bucataru and M.F. Dahl, \emph{Semi-basic 1-forms and {Helmholtz} conditions
  for the inverse problem of the calculus of variations}, Journal of Geometric
  Mechanics \textbf{1} (2009), no.~2, 159--180.

\bibitem[BD10a]{BucataruDahl:2008}
\bysame, \emph{A complete lift for semisprays}, International Journal of
  Geometric Methods in Modern Physics \textbf{7} (2010), no.~2, 267--287.

\bibitem[BD10b]{BucataruDahl:2008:conjugate}
\bysame, \emph{A geometric space without conjugate points}, Balkan Journal of
  Geometry and Its Applications \textbf{15} (2010), no.~1, 17--40.

\bibitem[Bes78]{Besse:1978}
A.L. Besse, \emph{Manifolds all of whose geodesics are closed}, Springer, 1978.

\bibitem[BM07]{BucataruMiron:2007}
I.~Bucataru and R.~Miron, \emph{{Finsler-Lagrange} geometry. {Applications} to
  dynamical systems}, Romanian {Academy}, 2007.

\bibitem[CP84]{CrampinPrince:1984}
M.~Crampin and G.E. Prince, \emph{The geodesic spray, the vertical projection,
  and {Raychaudhuri}'s equation}, General Relativity and Gravitation
  \textbf{16} (1984), no.~7, 675--689.

\bibitem[Dah08]{Dahl2007:Riccati}
M.F. Dahl, \emph{A geometric interpretation of the complex tensor {Riccati}
  equation for {Gaussian} beams}, Journal of non-linear mathematical physics
  \textbf{14} (2008), no.~1, 95--111.

\bibitem[Eis27]{Eisenhart}
L.P. Eisenhart, \emph{Non-riemannian geometry}, AMS Colloquium Publications,
  1927.

\bibitem[EJK98]{EhrlichJung:1998}
P.E. Ehrlich, Y.-T. Jung, and S.-B. Kim, \emph{Volume comparison theorems for
  {Lorentzian} manifolds}, Geometriae Dedicata \textbf{73} (1998), no.~1,
  39--56.

\bibitem[EK94]{EhrlichKim:1994}
P.E. Ehrlich and S.-B. Kim, \emph{From the {Riccati} inequality to the
  {Raychaudhuri} equation}, Contemporary Mathematics \textbf{170} (1994),
  65--78.

\bibitem[EO80]{EBurgSullivan}
J.-H. Eschenburg and J.J. O'Sullivan, \emph{Jacobi tensors and {Ricci}
  curvature}, Mathematische Annalen \textbf{252} (1980), no.~1, 1--26.

\bibitem[Hic65]{Hicks:1965}
N.J. Hicks, \emph{Notes on differential geometry}, Van Nostrand Reinhold
  Company, 1965.

\bibitem[JP00]{JeriePrince:2000}
M.~Jerie and G.E. Prince, \emph{A general {Raychaudhuri's} equation for
  second-order differential equations}, Journal of Geometry and Physics
  \textbf{34} (2000), 226--241.

\bibitem[JP02]{JeriePrince:2002}
\bysame, \emph{Jacobi fields and linear connections for arbitrary second order
  {ODE}'s}, Journal of Geometry and Physics \textbf{43} (2002), 351--370.

\bibitem[Kac05]{Kachalov:2005}
A.P. Kachalov, \emph{Gaussian beams, the {Hamilton-Jacobi} equations, and
  {Finsler} geometry}, Journal of Mathematical Sciences \textbf{127} (2005),
  no.~6, 2374--2388, English translation of article in Zapiski Nauchnykh
  Seminarov POMI, Vol. 297, 2003, pp. 66–92.

\bibitem[KKL01]{KKL:2001}
A.~Kachalov, Y.~Kurylev, and M.~Lassas, \emph{Inverse boundary spectral
  problems}, Chapman \& Hall/CRC, 2001.

\bibitem[KV86]{KowalskiVanhecke:1986}
O.~Kowalski and L.~Vanhecke, \emph{A new formula for the shape operator of a
  geodesic sphere and its applications}, Mathematische Zeitschrift \textbf{192}
  (1986), 613--625.

\bibitem[Lar96]{Larsen1996}
J.C. Larsen, \emph{The {Jacobi} map}, Journal of Geometry and Physics
  \textbf{30} (1996), 54--76.

\bibitem[Lew00]{Lewis:2001:GeomMax}
A.D. Lewis, \emph{The geometry of the maximum principle for affine connection
  control systems}, preprint (2000), 1--52.

\bibitem[MM63]{ManasseMisner:1963}
F.K. Manasse and C.W. Misner, \emph{Fermi normal coordinates and some basic
  concepts in differential geometry}, Journal of Mathematical Physics
  \textbf{4} (1963), no.~6, 735--745.

\bibitem[MV10]{MikesVanzurova:2010}
J.~Mike\v{s} and A.~Van\v{z}urov\'a, \emph{Reconstruction of the connection or
  metric from some partial information}, arXiv: 1006.3207 (2010), 1--9.

\bibitem[Sak96]{Sakai1992}
T.~Sakai, \emph{Riemannian geometry}, American Mathematical Society, 1996.

\bibitem[She01a]{Shen2001}
Z.~Shen, \emph{Differential geometry of spray and {Finsler} spaces}, Springer,
  2001.

\bibitem[She01b]{Shen2001b}
\bysame, \emph{Lectures of {Finsler} geometry}, World Scientific, 2001.

\bibitem[Spi79]{spivakI}
M.~Spivak, \emph{A comprehensive introduction to differential geometry, vol.
  1}, Publish or Perish, Inc., 1979.

\bibitem[YI73]{Yano1973}
K.~Yano and S.~Ishihara, \emph{Tangent and cotangent bundles}, Marcel Dekker
  Inc., 1973.

\end{thebibliography}

\providecommand{\bysame}{\leavevmode\hbox to3em{\hrulefill}\thinspace}
\providecommand{\MR}{\relax\ifhmode\unskip\space\fi MR }
\providecommand{\MRhref}[2]{%
  \href{http://www.ams.org/mathscinet-getitem?mr=#1}{#2}
}
\providecommand{\href}[2]{#2}

\end{document}